\documentclass[reqno]{amsart}

\usepackage{amsmath,amssymb,amsthm,amsfonts,
	enumerate,hyperref,cleveref}

\hypersetup{
	colorlinks=true,
	linkcolor=blue,
	filecolor=magenta,      
	urlcolor=cyan,
	pdftitle={Overleaf Example},
	pdfpagemode=FullScreen,
}
\usepackage{dsfont}
\usepackage[all]{xy}
\usepackage[T1]{fontenc}
\usepackage{tikz}
\usepackage{pgfplots}
\pgfplotsset{compat=1.15}
\usepackage{mathrsfs}
\usetikzlibrary{arrows}

\DeclareMathOperator{\im}{Im}

\theoremstyle{plain}
\newtheorem{theorem}{Theorem}[section]
\newtheorem{thrm}{Theorem}[section]

\newtheorem{prop}[thrm]{Proposition}
\newtheorem{lem}[thrm]{Lemma}

\theoremstyle{definition}

\crefrangeformat{equation}{#3(#1)#4--#5(#2)#6}

\crefname{thrm}{Theorem}{Theorems}
\crefname{theorem}{Theorem}{Theorems}
\crefname{lem}{Lemma}{Lemmas}
\crefname{cor}{Corollary}{Corollaries}
\crefname{prop}{Proposition}{Propositions}
\crefname{defn}{Definition}{Definitions}
\crefname{exm}{Example}{Examples}
\crefname{rem}{Remark}{Remarks}
\crefname{conj}{Conjecture}{Conjectures}
\crefname{quest}{Question}{Questions}
\crefname{section}{Section}{Sections}
\crefname{equation}{\unskip}{\unskip}
\crefname{enumi}{\unskip}{\unskip}
\crefname{subsection}{Subsection}{Subsections}

\newcommand{\re}{\Re\hbox{e}}

\begin{document}
\title[Additive mappings preserving orthogonality]{Additive mappings preserving orthogonality between complex inner product spaces}

\author[L. Li]{Lei Li}
\address[L. Li]{School of Mathematical Sciences and LPMC, Nankai University, 300071 Tianjin, China.}
\email{leilee@nankai.edu.cn}

\author[S. Liu]{Siyu Liu}
\address[S. Liu]{School of Mathematical Sciences and LPMC, Nankai University, 300071 Tianjin, China.}
\email{760659676@qq.com}

\author[A.M. Peralta]{Antonio M. Peralta}
\address[A.M. Peralta]{Instituto de Matem{\'a}ticas de la Universidad de Granada (IMAG), Departamento de An{\'a}lisis Matem{\'a}tico, Facultad de
	Ciencias, Universidad de Granada, 18071 Granada, Spain.}
\email{aperalta@ugr.es}

	
\begin{abstract} Let $H$ and $K$ be two complex inner product spaces with dim$(X)\geq 2$. We prove that for each non-zero additive mapping $A:H \to K$ with dense image the following statements are equivalent:  \begin{enumerate}[$(a)$] \item $A$ is (complex) linear or conjugate-linear mapping and there exists $\gamma >0$ such that $\| A (x) \| = \gamma  \|x\|$, for all $x\in X$, that is, $A$ is a positive scalar multiple of a linear or a conjugate-linear isometry;
\item There exists $\gamma_1 >0$ such that one of the next properties holds for all $x,y \in H$:
		\begin{enumerate}[$(b.1)$] \item $\langle A(x) |A(y)\rangle = \gamma_1  \langle x|y\rangle,$
			\item $\langle A(x) |A(y)\rangle = \gamma_1  \langle y|x \rangle;$
		\end{enumerate}
%\item $A$ is linear or conjugate-linear and preserves orthogonality in both directions; 
\item $A$ is linear or conjugate-linear and preserves orthogonality; 
\item $A$ is additive and preserves orthogonality in both directions; 
\item $A$ is additive and preserves orthogonality.
\end{enumerate}  This extends to the complex setting a recent generalization of the Koldobsky--Blanco--Turnšek theorem obtained by Wójcik for real normed spaces. 
\end{abstract}
	
	\maketitle
	
	%\tableofcontents
	
\section{Introduction} Elements $x,y$ in a real or complex inner product space $(H, \langle \cdot| \cdot \rangle)$ are called \emph{orthogonal} in the Euclidean sense ($x\perp_2 y$ in short) if $\langle x| y\rangle =0$. It is an intriguing question to determine how much information about $H$ is preserved by knowing all orthogonal pairs of elements in $H$. We are led to the following problem on preservers: If $K$ is another real or complex inner product space, a mapping $\Delta : H\to K$ \emph{preserves (Euclidean) orthogonality} if $\forall x,y \in  H, \  x \perp_2 y\Rightarrow  \Delta (x) \perp_2 \Delta (y).$ If the implication ``$\Rightarrow$'' is replaced with the equivalence ``$\Leftrightarrow$'' we say that $\Delta$ preserves (Euclidean) orthogonality in both directions. When can we conclude that such a mapping $\Delta$ is linear or conjugate linear isomorphism? Are the inner product spaces $X$ and $Y$ isometrically isomorphic? \smallskip

We should begin by observing that it is simply hopeless to characterize additive surjective maps preserving orthogonality in both directions on the one-dimensional complex Hilbert space $H = \mathbb{C}$, since every additive surjective mapping on $\mathbb{C}$ automatically preserves orthogonality in both directions. If we take a non-continuous bijective additive mapping $f:\mathbb{R}\to\mathbb{R}$, the natural extension $\tilde{f}: \mathbb{C}\to\mathbb{C},$ $\tilde{f} (\alpha + i \beta) = f(\alpha) + i f(\beta)$ is an additive bijection preserving orthogonality in both directions. We shall see next that this counterexample can only occur when the inner product space in the domain is one--dimensional.\smallskip

Concerning the above questions, a result by J. Chmieliński (see \cite[Theorem 1]{Chmielinski2005}) assures that for each non-zero mapping $T$ between two (real or complex) Hilbert spaces $H$ and $K$ the following statements are equivalent:\begin{enumerate}[$(a)$]
\item $T$ is linear and there exists $\gamma  >0$ such that $\|T (x)\| =\gamma \|x\|$ for all $x\in H$;
\item There exists $\gamma_1 >0$ such that $\langle T(x) |T(y)\rangle = \gamma_1  \langle x|y\rangle,$ for all $x,y \in H$;
\item  $T$ is linear and preserves orthogonality in both directions;
\item $T$ is linear and orthogonality preserving. 
\end{enumerate} Chmieliński also gave examples of non-linear (actually non-additive) maps on the $2$-dimensional real Hilbert space $\ell_2^2$.\smallskip

The conclusion in the just quoted result is closely related to the celebrated {Koldobsky--Blanco--Turnšek theorem}. It is known that elements $x,y$ in a Hilbert space $H$ are orthogonal if, and only if, they are \emph{Birkhoff or Birkhoff-James orthogonal} ($x\perp_{B} y$ in short), that is, for all $\alpha \in \mathbb{K}$ we have $\|x\| \leq  \|x  + \alpha y\|$. Birkhoff orthogonality makes sense for normed spaces. Let us observe that the underlying field $\mathbb{K}$ is crucial in the above definitions. Consider, for example, the complex Hilbert space $H=\ell^2_2$ with inner product $\langle (\lambda_1,\lambda_2)| (\mu_1,\mu_2)\rangle = \sum_{i} \lambda_i \overline{\mu_i}$, and the underlying real Hilbert space $H_{\mathbb{R}}$ with respect to the inner product $(\cdot|\cdot ) = \re \langle \cdot|\cdot\rangle$. It is easy to check that $((i,1) | (1,-i) ) =  \re (2 i ) =0,$ that is $(i,1) \perp_2 (1,-i)$ in $H_{\mathbb{R}},$ while $(i,1) \not\perp_2 (1,-i)$ in $H$. \smallskip

The celebrated \emph{Koldobsky--Blanco--Turnšek theorem} asserts that a non-zero linear mapping $T$ between two real or complex normed spaces $X$ and $Y$, preserves Birkhoff orthogonality if, and only if, there exists $\gamma >0$ such that $\| T (x) \| = \gamma  \|x\|$, for all $x\in X$ (see \cite{BlancoTurnsek} for the case of real normed spaces and \cite{Kol93} for complex normed spaces). In a more recent reference (cf. \cite[Theorem 3.1]{Wojcik}), P. Wójcik established a generalization of the {Koldobsky--Blanco--Turnšek theorem} by showing that for each non-zero additive mapping $A$ between two normed real spaces $X$ and $Y$ with dim$(X)\geq 2$, the following conditions are equivalent:  \begin{enumerate}[$(a)$]
	\item $A$ preserves Birkhoff orthogonality; 
	\item $A$ is a linear mapping and there exists $\gamma >0$ such that $\| T (x) \| = \gamma  \|x\|$, for all $x\in X$.
\end{enumerate} Obviously, Wójcik's theorem holds when $X$ and $Y$ are real inner product spaces, where it is more natural to speak about preservers of (Euclidean) orthogonality. However, the strong dependence on the base field of Birkhoff orthogonality and (Euclidean) orthogonality makes impossible to apply Wójcik's result for additive orthogonality preserver between complex Hilbert spaces (even more for additive preservers of Birkhoff orthogonality between complex normed spaces). It seems interesting to fill the existing gap and characterize all additive maps preserving orthogonality between complex Hilbert spaces. This short note is aimed to provide a full characterization of additive orthogonality preserving maps between complex inner product spaces, in this setting we show that, as it can be naturally expected, positive scalar multiples of  conjugate-linear isometries are also possible. In our main result we prove that if $A: H\to K$ is a non-zero additive mapping between complex inner product spaces, then $A$ is a positive scalar multiple of a (complex) linear or conjugate linear isometry if and only if it preserves (Euclidean) orthogonality (see \Cref{t Wojcik thm for complex Hilbert spaces}). We also prove that if $X$ and $Y$ are complex normed spaces with dim$(X)\geq 2$ and $X$ admits a conjugation $\tau$ (i.e., a period-$2$ conjugate linear isometry), then every additive mapping $A: X\to Y$ preserving Birkhoff orthogonality is real-linear. Furthermore, if $A$ is surjective and preserves Birkhoff orthogonality in both directions, then $A$ is a real-linear isomorphism and the underlying real normed spaces $X_{\mathbb{R}}$ and $Y_{\mathbb{R}}$ are isomorphic (see \Cref{p real-linearity additive op complex normed spaces}).\smallskip

\section{The results}

Before stating the desired characterization we include a technical lemma.

\begin{lem}\label{l complex multiples of a norm-one element} Let $A : H\to K$ be a real-linear mapping between two complex inner product spaces with dense image. Suppose additionally that $A$ preserves orthogonality. Then for each norm-one element $x_0\in H$ we have $A (\mathbb{C} x_0) \subseteq \mathbb{C} A(x_0)$.
\end{lem}

\begin{proof} We can suppose that $A(i x_0) \neq 0$, otherwise $A(\mathbb{C} x_0) = A(\mathbb{R} x_0) = \mathbb{R} A (x_0)\subseteq \mathbb{C} A (x_0)$.  \smallskip
	
Let us take $z\in K = \overline{A(H)}$ with $z \perp_2 A(x_0)$. By hypothesis, we can find a sequence $(x_n)$ in $H$ such that $A(x_n) \to z$. Since $H = \left(\mathbb{R} x_0 + \mathbb{R} (i x_0)\right) \oplus \{x_0\}^{\perp} = \left(\mathbb{R} x_0 + \mathbb{R} (i x_0)\right) \oplus \{i x_0\}^{\perp}$, there exist sequences $(\alpha_n)_n,$ $(\beta_n)_n$  in $\mathbb{R},$ and a sequence $(y_n)_n\subseteq \{i x_0\}^{\perp} = \{ x_0\}^{\perp}$ satisfying $\alpha_n x_0 + \beta_n i x_0 + y_n = x_n$ for all $n$, and hence $\alpha_n A (x_0)  + \beta_n  A(i x_0) + A(y_n) = A (\alpha_n x_0 + \beta_n i x_0 + y_n) \to z$.  
By hypotheses, $(\alpha_n A (x_0)  + \beta_n  A(i x_0)) \perp_2 A(y_n)$ for all natural $n$, and hence the sequences  $(\alpha_n A (x_0)  + \beta_n  A(i x_0))_n$ and $(A(y_n))_n$ must be bounded. \smallskip

If $A(x_0) =0 $, we have $A(H) = \mathbb{R} A(i x_0) \oplus^{\perp} A\left( \{x_0\}^{\perp} \right),$ and thus $$\begin{aligned}
K &= \overline{A(H)} = \overline{ \mathbb{R} A(i x_0) \oplus^{\perp} A\left( \{x_0\}^{\perp} \right) } = \mathbb{R} A(i x_0) \oplus^{\perp} \overline{A\left( \{x_0\}^{\perp} \right)} \\
&= \mathbb{R} A(i x_0) \oplus^{\perp} \overline{A\left( \{i x_0\}^{\perp} \right)} \subseteq \mathbb{R} A(i x_0) \oplus^{\perp}  \{A(i x_0)\}^{\perp} \subseteq K,
\end{aligned} $$ which is impossible. Therefore $A(x_0)\neq 0$. \smallskip

If $A(x_0)$ and $A(i x_0)$ are $\mathbb{R}$-linearly dependent we can write $A(i x_0) = t A(x_0)$ for some real $t$, and thus $A(i x_0 )\in \mathbb{C} A(x_0)$.\smallskip

We can now deal with the case that $A(x_0)$ and $A(i x_0)$ are $\mathbb{R}$-linearly independent. Since $(\alpha_n A (x_0)  + \beta_n  A(i x_0))_n$ is bounded, by basic theory of normed spaces the sequences, $(\alpha_n)_n$ and $(\beta_n)_n$ must be bounded. Up to taking appropriate subsequences, we can assume that $(\alpha_n)_n \to \alpha_0\in \mathbb{R}$ and $(\beta_n)_n\to \beta_0\in \mathbb{R}$, and hence $(A(y_n))_n\to z_1\in \overline{A(\{i x_0\}^{\perp})}\subseteq \{A(i x_0)\}^{\perp}\cap \{A(x_0)\}^{\perp}.$ We therefore have $z = \alpha_0 A(x_0) + \beta_0 A(i x_0) + z_1$.\smallskip 

We shall next show that $\beta_0 = 0$. Arguing by contradiction, we assume $\beta_0 \neq 0$. By applying that $z\perp_2 A(x_0)$, we deduce that \begin{equation}\label{eq to beta0 =0} \begin{aligned}
		0 = \langle z | A(x_0)\rangle &= \alpha_0 \|A(x_0)\|^2 +\beta_0  \langle A(i x_0) | A( x_0) \rangle+ \langle z_1 | A(x_0)\rangle \\
		&= \alpha_0 \|A(x_0)\|^2 +\beta_0  \langle A(i x_0) | A( x_0)\rangle,
	\end{aligned}
\end{equation} which implies that $\im \langle A(i x_0) | A( x_0)\rangle =0 $, equivalently, $\langle A(i x_0) | A( x_0)\rangle \in \mathbb{R}$.  This assures that $$A(i x_0) = s A(x_0) + z_2, \hbox{ with } s\in \mathbb{R},\hbox{ and } z_2\perp_2 A(x_0). $$ We therefore arrive to $$A(H) = (\mathbb{R} A( x_0)+ \mathbb{R} A(i x_0)) \oplus^{\perp} A\left( \{x_0\}^{\perp} \right) \subseteq \mathbb{R} A( x_0)\oplus^{\perp}  \{A(x_0)\}^{\perp} ,$$ and thus $K = \overline{A(H)} \subseteq \mathbb{R} A( x_0)\oplus^{\perp}  \{A(x_0)\}^{\perp} \subseteq K,$ which is also impossible, and hence $\beta_0 =0$, and by \eqref{eq to beta0 =0}, $\alpha_0 =0 $.\smallskip

All the previous conclusions lead to $z = z_1 \in \{A(i x_0)\}^{\perp}\cap \{A(x_0)\}^{\perp}$. We have then proved that $A(i x_0) \perp_2 z$ for every $z\perp_2 A(x_0)$, and then $A(i x_0 ) \in \mathbb{C} A(x_0)$. 
\end{proof}

We can now state the desired extension of the Blanco--Turnšek and Wójcik theorems.

\begin{theorem}\label{t Wojcik thm for complex Hilbert spaces} Let $A : H\to K$ be a non-zero mapping between two complex inner product spaces with dim$(H)\geq 2.$ Suppose that $A$ has dense image. Then the following statements are equivalent: \begin{enumerate}[$(a)$] \item $A$ is (complex) linear or conjugate-linear mapping and there exists $\gamma >0$ such that $\| A (x) \| = \gamma  \|x\|$, for all $x\in X$, that is, $A$ is a positive scalar multiple of a linear or a conjugate-linear isometry;
\item There exists $\gamma_1 >0$ such that one of the next properties holds for all $x,y \in H$:
\begin{enumerate}[$(b.1)$] \item $\langle A(x) |A(y)\rangle = \gamma_1  \langle x|y\rangle,$
	\item $\langle A(x) |A(y)\rangle = \gamma_1  \langle y|x \rangle;$
\end{enumerate}
\item $A$ is linear or conjugate-linear and preserves orthogonality in both directions; 
\item $A$ is linear or conjugate-linear and preserves orthogonality; 
\item $A$ is additive and preserves orthogonality in both directions; 
\item $A$ is additive and preserves orthogonality.
	\end{enumerate}
\end{theorem}

\begin{proof} The equivalence $(a)\Leftrightarrow (b)$ is known, actually $(b.1)$ (respectively, $(b.2)$) holds if, and only if, $A$ is linear (respectively, conjugate-linear). The implications $(b)\Rightarrow (c)\Rightarrow (d) \Rightarrow (f),$  $(c) \Rightarrow (e)$ and $(e) \Rightarrow (f)$ are clear.\smallskip
	
$(f) \Rightarrow (a)$ Suppose that $A$ is additive and preserves orthogonality. We shall first prove that $A$ is real-linear. One is first tempted to apply Wójcik's theorem \cite[Theorem 3.1]{Wojcik} to the additive mapping $A: H_{\mathbb{R}}\to K_{\mathbb{R}}$ regarded as a map between the underlying real inner product spaces $H_{\mathbb{R}}$ and $K_{\mathbb{R}}$ when equipped with the inner product $(x|y)= \re \langle a |b\rangle$. However, in such a case we can find points $x,y\in H$ with $(x|y)=0$ but $\langle x|y\rangle \neq 0$. We need a subtle argument to avoid the problem. \smallskip

Our goal is to show that $A (\alpha x_1) = \alpha A(x_1)$ for all $\alpha\in \mathbb{R}$ and $x_1\in H\backslash \{0\}$. By hypothesis, we can find $x_1\in H$ such that $A(x_1)\neq 0$. Since dim$(H)\geq 2$, there exists $x_2\in H\backslash \{0\}$ such that $x_1\perp_2 x_2$ in $H$. Let $H_1$ be the real subspace of $H$ generated by $\{x_1,x_2\}$, that is, $H_1 = \mathbb{R} x_1 \oplus \mathbb{R} x_2$. Since the elements in $H_1$ are real-linear combinations of $x_1$ and $x_2$, we can easily see that $\langle x| y\rangle =\re\langle x| y\rangle  = (x|y)$ for all $x,y\in H_1$, and thus elements in $H_1$ are orthogonal in $(H, \langle \cdot | \cdot\rangle)$ if, and only if, they are orthogonal in $(H_1, \re \langle \cdot| \cdot\rangle)$. We can therefore conclude that the mapping $A|_{H_1}: H_1 \to K_{\mathbb{R}}$ is additive and orthogonality preserving. So, by Wójcik's theorem \cite[Theorem 3.1]{Wojcik}, $A|_{H_1}$ is real-linear (and there exists a positive $\gamma$ satisfying $\|A(x) \| = \gamma \|x\|$ for all $x\in H_1$). We have also shown that $A(x_2)\neq 0$ for all non-zero $x_2$ with $\langle x_2 | x_1\rangle =0$. Replacing the roles of $x_1$ and $x_2$ we get $A(\lambda x_1) \neq 0$ for all $\lambda \in \mathbb{C}\backslash\{0\}$. If $x$ is a non-zero vector in $H$ we can find $\lambda\in \mathbb{C}$ and $x_2\perp_2 x_1$ such that $x = \lambda x_1 + x_2$ and $\|x\|^2 = \|\lambda x_1 \|^2 + \|x_2\|^2$, since $A(x)$ is the orthogonal sum of $A(\lambda x_1)$ and $A(x_2)$, it follows from the previous conclusions that $A(x)\neq 0$. Consequently, in the above arguments, $x_1$ can be replaced with any non-zero vector. The arbitrariness of $x_1$ allows us to conclude that $A$ is real-linear (and preserves orthogonality by hypothesis).\smallskip

We shall next prove that $A$ is (complex) linear or conjugate-linear. Observe first that, by \Cref{l complex multiples of a norm-one element},  $A(\mathbb{C} x_1) \subseteq \mathbb{C} A(x_1)$ for all norm-one element $x_1 \in H$. Fix a norm-one element $x_1\in H$. As before, since dim$(H)\geq 2$, we can find an orhtonormal system of the form $\{x_1, x_2\}$. A new application of \Cref{l complex multiples of a norm-one element} proves that $A (\mathbb{C} x_2) \subseteq  \mathbb{C} A ( x_2)$. Observe that, by hypotheses, each element in the set $\{A(x_1), A(i x_1)\}$ is orthogonal to every element in the set $\{A(x_2), A(i x_2)\}$.\smallskip

By applying that $A$ is real-linear we deduce that \begin{equation}\label{eq algebraic expression of A} A( \lambda_1 x_1 + \lambda_2 x_2  ) = \re (\lambda_1) A(x_1) + \im(\lambda_1) A(i x_1) + \re (\lambda_2) A(x_2) +  \im(\lambda_2) A(i x_2),
\end{equation} for all $\lambda_1,\lambda_2\in \mathbb{C}$. 
Observe that $i x_1 + i  x_2$ and $i x_1 -i  x_2$ are orthogonal vectors in $H$, and thus, by the hypothesis on $A$, we must have \begin{equation}\label{eq friday afternoon 1} \|A(i x_1)\|^2 - \|A(i x_2)\|^2   = \langle A( i x_1 + i  x_2 ) | A( i x_1 -i  x_2 ) \rangle=0.
\end{equation}

On the other hand, for every $\lambda_2 = \alpha + i \beta\in \mathbb{C}\backslash\{0\}$ and $\lambda_1 = s + i t\in \mathbb{C}$ the vectors $x_1 + \lambda_2 x_2$ and $\lambda_1 x_1 - \lambda_1 \overline{\lambda_2^{-1}} x_2$ are orthogonal in $H$, we can therefore conclude from the assumptions on $A$ that $A(x_1 + \lambda_2 x_2) \perp_2 A(\lambda_1 x_1 - \lambda_1 \overline{\lambda_2^{-1}} x_2)$, equivalently, $\langle A(x_1 + \lambda_2 x_2) | A(\lambda_1 x_1 - \lambda_1 \overline{\lambda_2^{-1}} x_2) \rangle=0,$ which expanded gives 
$$\begin{aligned} 0= &\re (\lambda_1) \|A(x_1)\|^2 + \im (\lambda_1) \langle A(x_1) | A(i x_1) \rangle +  \re (\lambda_2)  \re (-\lambda_1 \overline{\lambda_2^{-1}}) \|A(x_2)\|^2 \\
	&+  \re (\lambda_2)  \im (-\lambda_1 \overline{\lambda_2^{-1}}) \langle A(x_2) | A(i x_2) \rangle + \im (\lambda_2)  \re (-\lambda_1 \overline{\lambda_2^{-1}}) \langle A(i x_2) | A(x_2) \rangle \\
	&+ \im (\lambda_2) \im (-\lambda_1 \overline{\lambda_2^{-1}}) \| A(i x_2)\|^2.
\end{aligned}$$ If we rewrite the previous identity in terms of real and imaginary parts of $\lambda_2$ and $\lambda_1$ we arrive to 
\begin{equation}\label{eq key for real and imaginary parts} \begin{aligned} &0= s \|A(x_1)\|^2 + t \langle A(x_1) | A(i x_1) \rangle +  \alpha \frac{\beta t  -\alpha s}{\alpha^2+\beta^2}   \|A(x_2)\|^2 \\
		&-  \alpha  \frac{\beta s  +\alpha t}{\alpha^2+\beta^2}  \langle A(x_2) | A(i x_2) \rangle + \beta \frac{\beta t  -\alpha s}{\alpha^2+\beta^2}  \langle A(i x_2) | A(x_2) \rangle - \beta  \frac{\beta s  +\alpha t}{\alpha^2+\beta^2} \| A(i x_2)\|^2,
	\end{aligned} 
\end{equation} for all $\alpha,\beta,s,t\in \mathbb{R}$ with $\alpha^2+\beta^2 \neq 0$. Taking a quadruple of the form $(\alpha,\beta,s,t) = (1,0,-1,0)$ in \eqref{eq key for real and imaginary parts} we get \begin{equation}\label{eq friday afternoon 2} 0=- \|A(x_1)\|^2 + \|A(x_2)\|^2 .
\end{equation} By applying the latter equality to simplify \eqref{eq key for real and imaginary parts} we obtain 
\begin{equation}\label{equ 6 bis} \begin{aligned} &0=  \beta \frac{\beta s  + \alpha t}{\alpha^2+\beta^2}  \|A(x_1)\|^2 + t \langle A(x_1) | A(i x_1) \rangle-  \alpha  \frac{\beta s  +\alpha t}{\alpha^2+\beta^2}  \langle A(x_2) | A(i x_2) \rangle   \\
		&+ \beta \frac{\beta t  -\alpha s}{\alpha^2+\beta^2}  \langle A(i x_2) | A(x_2) \rangle - \beta  \frac{\beta s  +\alpha t}{\alpha^2+\beta^2} \| A(i x_2)\|^2,
	\end{aligned} 
\end{equation} for all $\alpha,\beta,s,t\in \mathbb{R}$ with $\alpha^2+\beta^2 \neq 0$. Take now $(\alpha,\beta,s,t) = (1,0,0,1)$ to obtain 
\begin{equation}\label{eq crossed inner products fisrt and second vectors} \langle A(x_1) | A(i x_1) \rangle -   \langle A(x_2) | A(i x_2) \rangle = 0.
\end{equation} By combining \eqref{eq crossed inner products fisrt and second vectors} and \eqref{equ 6 bis} we arrive to $$ \begin{aligned} &0=  \beta \frac{\beta s  + \alpha t}{\alpha^2+\beta^2}  \|A(x_1)\|^2 + \beta \frac{\beta t  -\alpha s}{\alpha^2+\beta^2} \langle A(x_1) | A(i x_1) \rangle  \\
	&+ \beta \frac{\beta t  -\alpha s}{\alpha^2+\beta^2}  \langle A(i x_2) | A(x_2) \rangle - \beta  \frac{\beta s  +\alpha t}{\alpha^2+\beta^2} \| A(i x_2)\|^2,
\end{aligned}$$ which leads to  
\begin{equation}\label{equ 6 bis bis} \begin{aligned} &0=  (\beta s  + \alpha t )  \left(\|A(x_1)\|^2 - \| A(i x_2)\|^2 \right) + 2 (\beta t  -\alpha s)  \re \langle A(x_1) | A(i x_1) \rangle,  
	\end{aligned} 
\end{equation} for all $\alpha,\beta,s,t\in \mathbb{R}$ with $\beta \neq 0$. In the cases $(\alpha,\beta,s,t) = (0,1,0,1)$ and $(\alpha,\beta,s,t) = (0,1,1,0)$ we get $$\re \langle A(x_1) | A(i x_1) \rangle =0, \hbox{ and } \|A(x_1)\|^2 =  \| A(i x_2)\|^2, \hbox{ respectively.}$$ Now, having in mind that $A(i x_1) \in \mathbb{C} A(x_1)$ and both are non-zero, $\|A(i x_1)\|= \|A(i x_2)\|$, and   $\|A(x_1)\| = \|A(x_2)\|$ (cf. \eqref{eq friday afternoon 1} and \eqref{eq friday afternoon 2}), it can be easily deduced that $A(i x_1)\in \{\pm i A(x_1)\}.$ The roles of $x_1$ and $x_2$ are clearly interchangeable, so $A(i x_2)\in \{\pm i A(x_2)\}.$ Let us write $A(i x_j) = \sigma_j (i) A(x_j)$ for $j\in\{1,2\}$ and $\sigma_j (i)\in\{\pm i\}$. Furthermore, it follows from the equality $\langle A(x_1) | A(i x_1) \rangle =   \langle A(x_2) | A(i x_2) \rangle $ (cf. \eqref{eq crossed inner products fisrt and second vectors}) that $\sigma_1 (i) = \sigma_2 (i)$. Back to \eqref{eq algebraic expression of A} we easily check that $A$ restricted to the linear complex span of $x_1$ and $x_2$ must be complex-linear or conjugate-linear.\smallskip

We have therefore shown that if $\{x_1,x_2\}$ is an orthonormal system in $H$, the restriction of $A$ to the complex-linear span of $\{x_1,x_2\}$ must be complex-linear or conjugate-linear.\smallskip

Let us pick a norm-one vector $z_0$ in $H$. It follows from the above that $A (i z_0) = i A(z_0)$ or $A (i z_0) = - i A(z_0)$. If the first case holds, for each norm one element $x$ in $H$, there exists another norm-one element $x_2$ in $H$ such that $\{z_0,x_2\}$ is an orthonormal system and $x$ belongs to its linear span. Since $A (i z_0) = i A(z_0)$, the previous conclusion proves that the restriction of $A$ to the linear span of $\{z_0,x_2\}$ is complex-linear, and thus $A(i x ) = i A(x)$, which implies that $A$ is complex-linear. If $A (i z_0) = - i A(z_0)$, similar arguments show that $A$ is conjugate-linear.\smallskip

We have implicitly shown that $A$ is a positive scalar multiple of an isometry. Namely, the conclusions in the first part of the proof guarantee that for each non-zero $x_1\in H$, there exists a positive $\gamma_1\in \mathbb{R}$ satisfying $\|A(x_1)\| = \gamma_1 \|x_1\|$. If $x_2$ is another non-zero vector in $H$ with $x_1\perp_2 x_2$ (we can clearly assume that $\|x_1 \|=\|x_2\|=1$), by considering the subspace $H_1 = \mathbb{R} x_1 \oplus \mathbb{R} x_2,$ as in the second part of the proof, we deduce from \eqref{eq friday afternoon 2} that $\gamma_1 = \gamma_1 \|x_1\| = \|A(x_1)\| = \|A(x_2)\| =  \gamma_2 \|x_2\|= \gamma_2$. We have therefore shown that if $\{x_1,x_2\}$ is an orthonormal system in $H$, we have $\|A(x_1)\| = \|A(x_2)\|$. Furthermore, since $A$ is linear or conjugate-linear, for each $x\in \mathbb{C} x_1 \oplus^{\ell_2} \mathbb{C} x_2$ we also have 
$$\begin{aligned} \|A(x)\|=&
\|A(\lambda_1 x_1 + \lambda_2 x_2 )\| \in \{ \|\lambda_1 A(x_1) + \lambda_2 A(x_2 )\|, \|\overline{\lambda_1} A(x_1) + \overline{\lambda_2} A(x_2 )\| \} \\
&= \left\{ \sqrt{|\lambda_1|^2 \|A(x_1)\|^2 + |\lambda_2|^2 \|A(x_2 )\|^2}, \sqrt{|\overline{\lambda_1}|^2 \|A(x_1)\|^2 + |\overline{\lambda_2}|^2 \|A(x_2 )\|^2} \right\} \\
&=  \|A(x_1)\| \sqrt{|\lambda_1|^2 + |\lambda_2|^2} =   \|A(x_1)\| \|\lambda_1 x_1 + \lambda_2 x_2 \| =  \|A(x_1)\| \|x\|.
\end{aligned} $$ Finally, a standard argument, like the one employed in the previous paragraph, assures that $\|A(x)\| = \|A(x_1)\|$ for every couple of norm-one vectors $x,x_1\in H$, which finishes the proof.\smallskip

Although the reasoning in the above paragraph is self-contained, and almost explicit from what we proved before, there is also another method to deduce that $A$ is a positive scalar multiple of a linear or conjugate-linear isometry. We already deduced that $A$ is linear or conjugate-linear (and preserves orthogonality). In the first case, we can apply Chmieliński' theorem \cite[Theorem 1]{Chmielinski2005} to $A$ and we get the desired conclusion. In the second case, let $\overline{K}$ be the complex inner product space obtained from $K$ by replacing the complex structure with the conjugate one, that is, $\lambda \odot x:= \overline{\lambda} x$ ($x\in \overline{K}$, $\lambda\in \mathbb{C}$) and inner product $\langle x|y\rangle^{rev} =\langle y|x\rangle$ ($x,y\in \overline{K}$). The mapping $A^{rev}: H\to \overline{K},$ $A^{rev} (x) = A(x)$ is a linear mapping preserving orthogonality, and hence Chmieliński' theorem proves that $A$ is positive scalar multiple of a conjugate-linear isometry.      
\end{proof}

The problem of determining when a real-linear mapping between complex Banach spaces is actually complex-linear or conjugate-linear is a topic studied in several contributions, for example Dang established in \cite[Proposition 2.6]{Dang92} that every real-linear surjective isometry between Cartan factors with rank $\geq 2$ must be either (complex) linear or conjugate-linear. In the case of rank-one Cartan factors (i.e. complex Hilbert spaces) the conclusion does not hold. For example, the mapping $R: \ell_2\to \ell_2$, $R ((\lambda_n)_n) = \left(\frac{1+(-1)^n}{2} \lambda_n + \frac{1-(-1)^n}{2} \overline{\lambda_n}\right)_n$ is a surjective real-linear which is not complex-linear nor conjugate-linear. In our result, the Hilbert spaces are rank-one Cartan factors and the mapping $A$ is not assumed to be surjective nor isometric, however, the hypothesis of being an orthogonality preserving additive mapping forces $A$ to be complex-linear or conjugate-linear. Similarly, the conclusions around Tingley's problem in the case of Hilbert spaces assert that every isometric mapping from the unit sphere of a Hilbert space $H$ ``into'' the unit sphere of another Hilbert space $K$ can be extended to a real-linear isometric mapping from $H$ into $K$ (see \cite[Theorem 2.2 and Corollary 2]{DingHilbert2002}), but nothing can be concluded about the (complex) linearity or conjugate-linearity of the extension.\smallskip

It is natural to ask whether a generalization of Wójcik theorem holds for complex normed spaces and Birkhoff orthogonality. Concerning this question, we can present some partial answer in the case that the domain is a complex normed space admitting a conjugation. We recall that a conjugation on a complex normed space $X$ is period-2 conjugate-linear linear $\tau: X\to X$. Define a conjugation $\tau^{\sharp} : X^*\to X^*$ given by $\tau^{\sharp} (\phi) (x) := \overline{\phi(\tau(x))}$ ($\phi\in X^*$, $x\in X$). The sets $X^{\tau} =\{x\in X: \tau (x) =x \}$ and $(X^*)^{\tau^{\sharp}} =\{\phi\in X^*: \tau^{\sharp} (\phi) =\phi \}$ are real-linear subspaces of $X$ and $X^*$, respectively. By construction $\phi (X^{\tau})\subseteq \mathbb{R}$ for all $\phi \in (X^{*})^{\tau^{\sharp}}$, and $X = X^{\tau}\oplus i X^{\tau}$. It is also known that the mapping $(X^*)^{\tau^{\sharp}}\ni \phi\mapsto \phi|_{X^{\tau}}=\re \phi|_{X^{\tau}}$ is a surjective linear isometry from $(X^*)^{\tau^{\sharp}}$ onto $(X^{\tau})^*$.%\smallskip 

\begin{prop}\label{p real-linearity additive op complex normed spaces} Let $X$ and $Y$ be complex normed spaces with dim$(X)\geq 2$, and assume that $X$ admits a conjugation $\tau$. Let $A: X\to Y$ be an additive mapping preserving Birkhoff orthogonality. Then $A$ is real-linear.  Furthermore, if $A$ is surjective and preserves Birkhoff orthogonality in both directions, then $A$ is a real-linear isomorphism and the underlying real normed spaces $X_{\mathbb{R}}$ and $Y_{\mathbb{R}}$ are isomorphic.   
\end{prop}

\begin{proof} A characterization of the Birkhoff orthogonality via Hahn-Banach theorem assures that elements $x,y$ in a real or complex normed space $Z$ are Birkhoff orthogonal if and only if there exists a norm-one functional $\phi\in Z^*$ satisfying $\phi (x) =\|x\|$ and $\phi (y)  =0$. Consider the real subspace $X^{\tau}$. If $x\perp_{B} y$ in $X^{\tau}$ there exists a norm-one functional $\phi\in (X^{\tau})^*\equiv (X^*)^{\tau^{\sharp}}$ satisfying $\phi (x) =\|x\|$ and $\phi (y) =\overline{\phi (\tau(y))} = \re \phi (y) =0$. In particular, $\phi$ is a norm-one functional in $X^{*}$ (with $\tau^{\sharp} (\phi) =\phi$), and hence $x\perp_{B} y$ in $X$. Therefore, $x\perp_{B} y$ in $X^{\tau}$ if, and only if, $x\perp_{B} y$ in $X$. This implies that $A|_{X^{\tau}} : X^{\tau}\to Y_{\mathbb{R}}$ is an additive mapping preserving Birkhoff orthogonality between two normed real spaces. So, Wójcik's theorem \cite[Theorem 3.1]{Wojcik} implies that $A|_{X^{\tau}}$ is real-linear and there exists a positive constant $\gamma_1$ satisfying $\| A(x)\| = \gamma_1 \|x\|$ for all $x\in X^{\tau}$. 
By applying a similar argument to the real subspace $i X^{\tau}$, whose dual space can be identified with $i (X^*)^{\tau^{\sharp}} = \{\phi\in X^* : \tau^{\sharp} (\phi) =-\phi\}$ (or by simply replacing $\tau$ with $\tau_1 =- \tau,$ and apply the above argument to $X^{\tau_1}= i X^{\tau}$), we deduce via Wójcik's theorem \cite[Theorem 3.1]{Wojcik} that $A|_{i X^{\tau}}$ is real-linear too, and there exists a positive constant $\gamma_2$ satisfying $\| A( i y)\| = \gamma_2 \|i y\| = \|y\|$ for all $y\in X^{\tau}$. 
Finally, since $X = X^{\tau}\oplus i X^{\tau}$, the mapping $A$ must be real linear. Furthermore, $$\|A(x + i y)\| \leq \|A(x)\| + \|A(i y)\| =  \gamma_1 \|x\| + \gamma_2 \|y\| \leq (\gamma_1 + \gamma_2) \|x+ i y \|,$$ for all $x+i y \in  X^{\tau}\oplus i X^{\tau}= X$. \smallskip

Assume now that $A$ is surjective and preserves Birkhoff orthogonality in both directions (and, of course, real-linear). We first observe that $A$ is injective. Namely, if $A(x)=0$, it follows that $A(x)\perp_{B} A(z)$ for all $z\in X$, and hence $x\perp_{B} z$ for all $z\in X$, which clearly gives $x=0.$ Therefore, $A$ is a continuous real linear bijection. By replacing $A$ with $A^{-1}$ we deduce that $A$ is a real linear isomorphism of normed spaces.    
\end{proof}

Paraphrasing R. Tanaka \cite{Tanaka22Indag}, the conclusion in our \Cref{t Wojcik thm for complex Hilbert spaces} assures that the vector addition and the relationship of orthogonality determine the entire structure of a complex inner product space, while by \Cref{p real-linearity additive op complex normed spaces} vector addition and the relationship of Birkhoff orthogonality determine isomorphically the entire structure of a complex normed space admitting a conjugation. It can added that in \cite{Tanaka22}, Tanaka studies the problem whether the existence of a (possibly non-additive) bijection $\Delta$ preserving Birkhoff–James orthogonality in both directions between two real Banach spaces $X$ and $Y$ assures the existence of a linear isomorphism $\Phi$ between $X$ and $Y$, however no conclusion is obtained on the mapping $\Delta$ itself. Positive answers are known when $X$ is finite-dimensional, or when  $X$ and $Y$ are reflexive and smooth, or when $X$ is a Hilbert space with dim$(X)\geq 3$ (cf. \cite{Tanaka22}). \smallskip \smallskip

\textbf{Acknowledgements} L. Li was supported by National Natural Science Foundation of China (Grant No. 12171251). A.M. Peralta supported by grant PID2021-122126NB-C31 funded by MICIU/AEI/10.13039/501100011033 and by ERDF/EU, by Junta de Andalucía grant FQM375, IMAG--Mar{\'i}a de Maeztu grant CEX2020-001105-M/AEI/10.13039/501100011033 and (MOST) Ministry of Science and Technology of China grant G2023125007L. \smallskip

Part of this work was completed during a visit of A.M. Peralta to Nankai University and the Chern Institute of Mathematics, which he thanks for the hospitality.\medskip

We thank B. Kuzma for addressing our attention to Wójcik's paper \cite{Wojcik}.\smallskip

\subsection*{Statements and Declarations} 

All authors declare that they have no conflicts of interest to disclose.

\subsection*{Data availability}

There is no data associated with this article.

%\bibliographystyle{plain}
%\bibliography{bibl042024.bib}

\begin{thebibliography}{10}
	
%	\bibitem{Alfsen_Shultz_Stomer_AdvMath_1978}	E.M. Alfsen, F.W. Shultz, and E.~St{\o}rmer.	\newblock A {G}elfand-{N}eumark theorem for {J}ordan algebras.	\newblock {\em Advances in Math.}, \textbf{28}(1):11--56, 1978.
	
%	\bibitem{AyuArzi2016Rickart}	S.A. Ayupov and F.N. Arzikulov.	\newblock {Jordan counterparts of Rickart and Baer $^*$-algebras}.\newblock {\em Uzbek. Mat. Zh.}, \textbf{13}(1):220--233, 2016.
	
\bibitem{BlancoTurnsek} A. Blanco, A. Turnšek, On maps that preserve orthogonality in normed spaces, \emph{Proc. Roy. Soc. Edinburgh Sect. A} \textbf{136} (2006), 709--716. 
	
%\bibitem{Bar_Tim_MathScand_1986}	T.~Barton and Richard~M. Timoney.	\newblock Weak{$^\ast$}-continuity of {J}ordan triple products and its applications. \newblock {\em Math. Scand.}, \textbf{59}(2):177--191, 1986.
	
%\bibitem{Battaglia_1991}	M.~Battaglia. \newblock Order theoretic type decomposition of {JBW}-triples.	\newblock {\em Q. J. Math., Oxf. II. Ser.}, \textbf{42}(166):129--147, 1991.
	
%	\bibitem{BeLoPeRo} J. Becerra Guerrero, G. L{\'o}pez P{\'e}rez, A. M. Peralta, A. Rodr{\'\i}guez-Palacios, A. Relatively weakly open sets in closed balls of Banach spaces, and real ${\rm JB}^*$-triples of finite rank, \emph{Math. Ann.} \textbf{330},  no. 1, 45--58 (2004).
	
%	\bibitem{BraunKaupUpmeier78}	R.~Braun, W.~Kaup, and H.~Upmeier.	\newblock A holomorphic characterization of {Jordan} {{\(C^*\)}}-algebras.	\newblock {\em Math. Z.}, \textbf{161}:277--290, 1978.
	
%\bibitem{BuChu92} L.J. Bunce, Ch.-H. Chu,  Compact  operations, multipliers and Radon-Nikodym property in JB$^*$-triples, \emph{Pacific J. Math.} \textbf{153}, 249--265 (1992).
	
%\bibitem{BunFerMartPe06} L.J. Bunce, F.J. Fernández-Polo, J.M. Moreno, A.M. Peralta, A Saitô–Tomita–Lusin theorem for JB$^*$-triples and	applications, \emph{Quart. J. Math. Oxford Ser.} \textbf{57} (2006), 37-48.
	
%\bibitem{BFPGMP08}	M.~Burgos, F.J. Fern{\'a}ndez-Polo, J.J. Garc{\'e}s, J.M. Mart{\'{\i}}nez, and A.M. Peralta. \newblock {Orthogonality preservers in C\(^*\)-algebras, JB\(^*\)-algebras and JB\(^*\)-triples}. \newblock {\em J. Math. Anal. Appl.}, \textbf{348}(1):220--233, 2008.
	
%\bibitem{BurKaMoPeRa} M.~Burgos, A.~Kaidi, A.~Morales, A.M. Peralta, and M.I. Ram{\'\i}rez. \newblock {von Neumann regularity and quadratic conorms in JB$^*$-triples and	C$^*$-algebras}. \newblock {\em Acta Math. Sinica (English Series)}, \textbf{24}(2):185--200, 2008.
	
%\bibitem{CabreraPalaciosBook}	M. Cabrera Garc{\'{\i}}a, A. Rodríguez Palacios,	\newblock {\em Non-associative normed algebras. {Volume} 1. {The} {Vidav}-{Palmer} and {Gelfand}-{Naimark} theorems}, volume \textbf{154} of {\em Encycl. Math. Appl.} \newblock Cambridge: Cambridge University Press, 2014.
	
%\bibitem{Chabb2017} F.~Chabbabi. \newblock Product commuting maps with the {{\(\lambda\)}}-{Aluthge} transform. \newblock {\em J. Math. Anal. Appl.}, \textbf{449}(1):589--600, 2017.
	
%\bibitem{ChabbMbekhta2017}	F.~Chabbabi and M.~Mbekhta. \newblock Jordan product maps commuting with the {{\(\lambda\)}}-{Aluthge} transform. \newblock {\em J. Math. Anal. Appl.}, \textbf{450}(1):293--313, 2017.
	
%\bibitem{ChabbMbekhta2022corrigendum}	F.~Chabbabi and M.~Mbekhta.	\newblock Corrigendum to: ``{Jordan} product maps commuting with the {{\(\lambda \)}}-{Aluthge} transform''.	\newblock {\em J. Math. Anal. Appl.}, \textbf{511}(1):1, 2022. \newblock Id/No 126029.
	
\bibitem{Chmielinski2005} J. Chmieliński, Linear mappings approximately preserving orthogonality, \emph{J. Math. Anal. Appl.} \textbf{304}  (2005), No. 1, 158--169. 
	
%	\bibitem{CUetoPeralta2022LMA}	M.~Cueto-Avellaneda and A.M. Peralta.\newblock Can one identify two unital {{\(\mathrm{JB}^*\)}}-algebras by the metric spaces determined by their sets of unitaries?	\newblock {\em Linear Multilinear Algebra}, \textbf{70}(22):7702--7727, 2022.
	
\bibitem{Dang92} T. Dang, Real isometries between JB$^*$-triples, \emph{Proc. Am. Math. Soc.} \textbf{114} (1992), No. 4, 971--980. 
	
%	\bibitem{Dang_Friedm_MathScand_1987} T.~Dang and Y.~Friedman.	\newblock Classification of {{\(JBW^*\)}}-triple factors and applications. \newblock {\em Math. Scand.}, \textbf{61}(2):292--330, 1987.
	
%	\bibitem{Dineen_theseconddual_1986} 	S.~Dineen. \newblock The second dual of a {JB{$^*$}} triple system. \newblock In {\em Complex analysis, functional analysis and approximation theory ({C}ampinas, 1984)}, volume \textbf{125} of {\em North-Holland Math. Stud.}, pages 67--69. North-Holland, Amsterdam, 1986.
		
\bibitem{DingHilbert2002} G. Ding, The $1$-Lipschitz mapping between the unit spheres of two Hilbert spaces can be extended to a real-linear isometry of the whole space, \emph{Sci. China, Ser. A} \textbf{45} (2002), No. 4, 479--483.
	
%\bibitem{DymMcKeanBook1972}	H.~Dym and H.P. McKean. \newblock Fourier series and integrals. \newblock Probability and {Mathematical} {Statistics}. {Vol}. 14. {New}{York}-{London}: {Academic} {Press}, 1972.
	
	
%\bibitem{EDW_RUTT_JLMS_1988}	C.M. Edwards and G.T. R\"{u}ttimann.	\newblock On the facial structure of the unit balls in a {${\rm JBW}^*$}-triple and its predual. \newblock {\em J. London Math. Soc. (2)}, \textbf{38}(2):317--332, 1988.
	
%\bibitem{EssPe2018} A.B.A. Essaleh and A.M. Peralta. \newblock Preservers of {{\(\lambda\)}}-{Aluthge} transforms.	\newblock {\em Linear Algebra Appl.}, \textbf{554}:86--119, 2018.
	
%\bibitem{FerGarSanSi92} A. Fernández L{\'o}pez, E. García Rus, E. Sánchez Campos, M. Siles Molina. \newblock {Strong regularity and generalized inverses in Jordan systems}. \newblock {\em Comm. Algebra}, \textbf{20}(7):1917--1936, 1992.	
	
%\bibitem{FerPe2007}	F.J. Fernández-Polo, A.M. Peralta, Compact tripotents and the Stone–Weierstrass Theorem for C$^*$-algebras and JB$^*$-triples, \emph{J. Operator Theor.} \textbf{58} (1) (2007), 157-173.

%\bibitem{FerPe2010} F.J. Fern{\'a}ndez-Polo and A.M. Peralta, Non-commutative generalisations of Urysohn's lemma and hereditary inner ideals, \emph{Journal of Functional Analysis} \textbf{259} (2010), 343-358
	
%\bibitem{Polo_Peralta_AdvMath_2018}	F.J. Fern{\'a}ndez-Polo and A.M. Peralta.	\newblock Low rank compact operators and {Tingley}'s problem, \newblock {\em Adv. Math.}, \textbf{338}:1--40, 2018.
	
%\bibitem{FriedmanHakeda1988}	Y.~Friedman and J.~Hakeda. \newblock Additivity of quadratic maps. \newblock {\em Publ. Res. Inst. Math. Sci.}, \textbf{24}(5):707--722, 1988.
	
%\bibitem{FriHak1988}	Y.~Friedman and J.~Hakeda. \newblock Additivity of quadratic maps.	\newblock {\em Publ. Res. Inst. Math. Sci.}, \textbf{24}(5):707--722, 1988.
	
%\bibitem{Fried_Peralta_Ann_Math_Phys_2022}	Y.~Friedman and A.M. Peralta. \newblock Representation of symmetry transformations on the sets of tripotents of spin and {Cartan} factors.	\newblock {\em Anal. Math. Phys.}, \textbf{12}(1):52, 2022. \newblock Id/No 37.

%\bibitem{FriedmanRusso82TAMS} Y. Friedman, B. Russo, Contractive projections on $C_0(K)$, \emph{Trans. Am. Math. Soc.}  (1982), \textbf{273}, 57-73.
	
%\bibitem{FriRuCommutative}  Y. Friedman, B. Russo, Function representation of commutative operator triple systems, \emph{J. Lond. Math. Soc.}, II. Ser.  (1983), \textbf{27}, 513-524.
	
%\bibitem{FriedRusso_Predual}	Y.~Friedman and B.~Russo. \newblock Structure of the predual of a {JBW\(^*\)}-triple. \newblock {\em J. Reine Angew. Math.} (1985), \textbf{356}:67--89.
	
%\bibitem{FriedRusso_GelfNaim}	Y.~Friedman and B.~Russo. \newblock The {Gelfand}-{Naimark} theorem for {JB\(^*\)}-triples. \newblock {\em Duke Math. J.}, \textbf{53}:139--148, 1986.

	
%\bibitem{GarLiPeSu} J.J. Garcés, L. Li, A.M. Peralta, S. Su, Maps preserving the truncation of triple products on Cartan factors. Preprint 2024,  arXiv:2405.13489. 
	
%\bibitem{GudderStrawther75} S. Gudder, D. Strawther, Orthogonally additive and orthogonally increasing functions on vector spaces, \emph{Pacific J. Math.} \textbf{58} (1975), 427--436.
	
%\bibitem{Hakeda1986} J.~Hakeda. \newblock Additivity of *-semigroup isomorphisms among *-algebras. \newblock {\em Bull. Lond. Math. Soc.}, \textbf{18}:51--56, 1986.
	
%\bibitem{HakedaSaito1986} J.~Hakeda and K.~Sait{\^o}. \newblock Additivity of {Jordan} *-maps between operator algebras. \newblock {\em J. Math. Soc. Japan}, \textbf{38}:403--408, 1986.
	
%\bibitem{Hakeda1986Jordan}	J\^{o}suke Hakeda. \newblock Additivity of {J}ordan{$^\ast$}-maps on {$AW^\ast$}-algebras. 	\newblock {\em Proc. Amer. Math. Soc.}, 96(3):413--420, 1986.
	
%\bibitem{HamKalPe2023} J.~Hamhalter, O.F.K. Kalenda, and A.M. Peralta.	\newblock Determinants in {Jordan} matrix algebras. \newblock {\em Linear Multilinear Algebra}, 71(6):961--1002, 2023.
	
%\bibitem{Horn_MathScand_1987} G.~Horn. \newblock Characterization of the predual and ideal structure of a {{\(JBW^*\)}}- triple.	\newblock {\em Math. Scand.}, \textbf{61}(1):117--133, 1987.
	
%\bibitem{Jia_Shi_Ji_AnnFunctAnn_2022}	X.~Jia, W.~Shi, and G.~Ji.	\newblock Maps preserving the truncation of products of operators.	\newblock {\em Ann. Funct. Anal.}, \textbf{13}(3):17, 2022. \newblock Id/No 40.

%\bibitem{Kal_Peralta_Ann_Math_Phys_2021}	O.F.K. Kalenda and A.M. Peralta.	\newblock Extension of isometries from the unit sphere of a rank-2 {Cartan}	factor. \newblock {\em Anal. Math. Phys.}, \textbf{11}(1):25, 2021.	\newblock Id/No 15.
	
%	\bibitem{Kaup77MathAnn} W. Kaup, Algebraic characterization of symmetric complex Banach manifolds, \emph{Math. Ann.}  (1977), \textbf{228}, 39-64.
	
%	\bibitem{Kaup_RiemanMap}	W.~Kaup. \newblock A {R}iemann mapping theorem for bounded symmetric domains. \newblock {\em Math Z.}, \textbf{183}:503--530, 1983.
	
%\bibitem{Kaup_Sing_Val}	W.~Kaup. \newblock On spectral and singular values in {JB\(^*\)}-triples.	\newblock {\em Proc. R. Ir. Acad., Sect. A}, \textbf{96}(1):95--103, 1996.
	
%	\bibitem{Ka96}	W.~Kaup. \newblock {On spectral and singular values in JB$^*$-triples}. \newblock {\em Proc. Roy. Irish Acad. Sect. A}, \textbf{96}(1):95--103, 1996.
	
%	\bibitem{Kaup_ManMath_1997} W.~Kaup. \newblock On real {Cartan} factors. \newblock {\em Manuscr. Math.}, \textbf{92}(2):191--222, 1997.
	
\bibitem{Kol93} A. Koldobsky, Operators preserving orthogonality are isometries, \emph{Proc. Roy. Soc. Edinburgh Sect. A} \textbf{123} (1993), 835--837.

	
%\bibitem{FerGarSanSi94}	A.F. L{\'o}pez, E.G. Rus, E.S. Campos, and M.S. Molina. \newblock {Strong regularity and Hermitian Hilbert triples}. \newblock {\em Quart. J. Math. Oxford Ser.(2)}, \textbf{45}(177):43--55, 1994.
	
%\bibitem{MaoJi2024}	Y.~Mao and G.~Ji. \newblock {Truncations of operators in {$\mathcal{B}(H)$} and their preservers}. \newblock {\em Adv. Operator Theory}, \textbf{9}:13, 2024. \newblock Id/No 32.
	
%\bibitem{OmlaSemrl93} M. Omladi\v{c}, P. \v{S}emrl, Additive mappings preserving operators of rank one, \emph{Linear Algebra Appl.} (1993), \textbf{182}, 239-256.
	
%bibitem{Per_RIMS_2023} A.M. Peralta.	\newblock Maps preserving triple transition pseudo-probabilities. \newblock {\em RIMS K{\^o}ky{\^u}roku Bessatsu}, \textbf{B93}:1--28, 2023.
	
%\bibitem{Peralta_ResMath_2023}	A.M. Peralta.	\newblock Preservers of triple transition pseudo-probabilities in connection	with orthogonality preservers and surjective isometries.	\newblock {\em Result. Math.}, \textbf{78}(2):23, 2023. \newblock Id/No 51.
	
%\bibitem{PeRu2014}	A.M. Peralta and B.~Russo. \newblock {Automatic continuity of triple derivations on C$^*$-algebras and		JB$^*$-triples}. \newblock {\em J. Algebra}, \textbf{399}:960--977, 2014.

\bibitem{Ratz} J. Rätz, On orthogonally additive mappings, \emph{Aequationes Math.} \textbf{28} (1985), 35--49.

\bibitem{Tanaka22} R. Tanaka, Nonlinear equivalence of Banach spaces based on Birkhoff-James orthogonality, \emph{J. Math. Anal. Appl.} \textbf{505} (2022), No. 1, Article ID 125444, 12 p.

\bibitem{Tanaka22Indag} R. Tanaka, On Birkhoff-James orthogonality preservers between real non-isometric Banach spaces, \emph{Indag. Math., New Ser.}  \textbf{33} (2022), No. 6, 1125-1136.

\bibitem{Wojcik} P. Wójcik, Mappings preserving B-orthogonality, \emph{Indag. Math., New Ser.} \textbf{30}  (2019), No. 1, 197--200. 
	
%	\bibitem{Wright1977}	J.D.M. Wright. \newblock Jordan {{\(C^*\)}}-algebras.	\newblock {\em Mich. Math. J.}, \textbf{24}:291--302, 1977.
	
%\bibitem{Yao_Ji_JMathResAppl_2022} 	J.~Yao and G.~Ji. \newblock Additive maps preserving the truncation of operators. \newblock {\em J. Math. Res. Appl.}, \textbf{42}(1):89--94, 2022.
	
\end{thebibliography}

%\bibliographystyle{acm}

\end{document}